\magnification=\magstep1
\input amstex
\documentstyle{amsppt}
\pageheight{8.5truein}
\NoBlackBoxes
\leftheadtext{Geometric Characterization of Solitons}
\rightheadtext{Soffer}
\def\varep{\varepsilon}
\topmatter
\title  Geometric Characterization of Solitons
\endtitle
\author  A. Soffer\endauthor
\address Mathematics Department,
 Rutgers University, New Brunswick, NJ 08903
\endaddress
\email soffer\@math.rutgers.edu
\endemail
\abstract I show that $H^1$ solutions of the nonlinear Schroedinger
equation which are incoming converge to a soliton, in the radial
case.
\endabstract
\endtopmatter

\document

\head Section 1\endhead
\subhead
1. Introduction\endsubhead
\medskip

The aim of this work is to describe a geometric definition of
localized solutions of NLS.

In the linear case we have the RAGE Theorem, which relates localized
solutions to the pure point spectrum of the Hamiltonian: Localized
solutions of the linear Schroedinger equation are linear
combinations of $L^2$ eigenfunctions of the Hamiltonian. In
particular, they are almost periodic functions of time. For the
nonlinear case see [Sig, Sof, Tao].

 The question arises as to what
is the analog of the bound states of a linear equation, in the
nonlinear case.

Here, I will show that solitons appear naturally from geometric
considerations. It lends support to the conjecture that all generic
outgoing states of NLS are solitons and free waves. That is, I will
show that if the solution of NLS is purely incoming, up to $L^1(dt)$
corrections, then the solution converges to a soliton, in any
compact region around the origin.

The method of proof is based on and motivated by the hydrodynamic
reformulation of the Schroedinger equation.

 The incoming wave
condition is then written in terms of the notion of flux through
surfaces around the origin.

It is then shown how to rigorously use the hydrodynamic formulation,
by restricting the analysis to topologically trivial domains of
space-time where the solution is nonvanishing.

 The solution in such
regions can then be uniquely written in the polar form, with
continuous phase function. This, together with the a-priori $H^1$
bound is then used to construct the velocity function, and make
sense of the related Euler type equation it satisfies. It should be
noted that here we do not use the semiclassical limit, so the
"Quantum Potential" term is not ignored.

Consider the NLS in 1 or 3 dimensions (for simplicity):

$$
i\frac {\partial \psi}{\partial t} = - \Delta \psi + F(|\psi|) \psi
\tag 1.1
$$
$\psi(t=0) = \psi_0 \in H^1_{\text{radial}}\cap L^2(\Bbb R^n, |x|
d^nx)$
that is
$$
\| |x|\psi_0\|_{L^2} + \|\nabla \psi_0\|_{L^2} + \|\psi_0\|_{L^2} <
\infty.\tag 1.2
$$
The nonlinearity $F(|\psi|)$ is assumed to be of the RSS type [RSS]
which guaranties global existence and stable soliton (ground state)
solutions, again for simplicity. The property of stability is not
used.

We will now make the following main time-dependent a-priori
assumptions on the solutions of (1) with $\psi_0$ as initial data:

\subhead $H^1$ (Energy) Boundedness\endsubhead

$$
\| \psi(t)\|_{H^1} \leq E < \infty \text{ uniformly in } t.\tag 1.3
$$

This assumption is of course verified whenever global existence is
proved.

\subhead Incoming Wave Condition (IWC)\endsubhead

$$
i \hat r \cdot (\psi\nabla\bar\psi - \bar \psi \nabla\psi) \leq 0
\text{ for all times, any } r.\tag 1.4
$$

Here $\hat r = \overarrow r/r = x/|x|$. The above inequality is to
be understood as holding almost everywhere in space. That this
notion makes sense follows since the above expression is an $L^1$
function by (1.3).

 This condition will be relaxed in various ways
later. We now combine the above two conditions with the energy and
dilation identities:

$$
\frac{d}{dt}  \left\{ \frac{1}{2} \int_{\Bbb R^n} |\nabla \psi |^2 + \int
\tilde F (|\psi|)|\psi|^2\right\} = 0 \tag 1.5a
$$
$$
\frac{d}{dt}  \|\psi\|_{L^2(\Bbb R^n)}^2 = 0  \tag 1.5b
$$
$$
\frac{d}{dt}  \| |x|\psi\|^2_{L^2(\Bbb R^n)} = (\psi, i[-\Delta,
|x|^2]\psi)
   = 2(\psi, (-i\nabla_x\cdot x - i x \cdot \nabla_x)\psi) = 4(\psi,
  A\psi) \leq 0. \tag 1.5c
$$
$$
\text{where } A\equiv (-ix\cdot \nabla_x - i \nabla_x \cdot x) / 2
   \tag 1.6
$$
where we used the IWC at the last step.

$\tilde F(|\psi|)$ is obtained from $F(|\psi|)$ in a known way.
>From the last identity we get
$$
\sup_t\int |x\psi|^2 d^n x \leq \int |x\psi_0|^2 d^n x < \infty. \tag 1.7
$$
Since we also have $\psi \in H^1, \sup_t \|\psi\|_{H^1} < E \text{
it follows }$

\proclaim{Proposition 1.1}

a) The trajectory $\{ \psi (t)\}^\infty_{t=0}$ is precompact in $H^s$,
for all $ 0\leq s < 1$.

b) Define $|\psi|^2 = \rho \text{ and } |\psi| = \sqrt{ \rho} = \eta$
then, since
$$
\big|\nabla |\psi|\big|\leq |\nabla \psi| \text{ a.e.,}
$$
$$
\{ |\psi (t)| = \eta (t)\}^\infty_{t=0}\text{ is precompact in }H^s,
0\leq s < 1.
$$
\endproclaim

\noindent{\bf Remarks}

$\bullet$ We will repeatedly use the fact that weak convergence of
$\eta_t$ implies strong convergence in $H^s, 0\leq s < 1$, and in
$L^p, 2 \leq p < \frac{2n}{n-2}$, due to the uniform $H^1$ bound.

$\bullet$ The above arguments extend to any monotonic increasing
function of $r$ replacing $r$ in the dilation identity, with
appropriate assumption on the localization of $\psi_0$:

$$
\| |x|^\delta \psi_0 \| < \infty. \qquad   0< \delta .
$$

\head Section 2 - Convergence of $\rho_t, \eta_t$\endhead

\subhead The vector field $\underline{v}$\endsubhead

For $\psi$ as above consider the current density $\vec J$.
$$
\vec J = - \frac{i}{2} (\bar \psi \nabla \psi - \psi \nabla \bar
\psi).\tag 2.1
$$
Since $\psi, \nabla \psi$ are in $L^2$ by our $H^1$ bound, $\vec J$ is
in $L^1$, uniformly in $t$.  For $(x, t)\in \Bbb R^n\times \Bbb R$
such that $\psi(x, t)\neq 0,$ we can define
$$
|\psi |^{-2}\vec J = \rho^{-1} \vec J = \tilde{\underline{v}}(x,
 t).\tag 2.2
$$
$\tilde{\underline{v}}(x,t)$ is defined on the set $Q=\{(x, t)|\psi(x,
t) \neq 0\}$.  We also denote, for each $t$, by $N_t$ the set where
$\psi(x, t)=0$.  In one dimension $\psi$ is bounded and continuous.
Moreover, $Q$ is open, $N^c_t$ is open.  So, for each point $(x, t)\in
Q$, we can find a neighborhood $U$ such that $|\psi(x', t')|>\delta,
 (x', t')\in U$, by (joint) continuity, see below.

Since $\vec J\in L^2(\Bbb R)$, it follows that $|\psi|^{-2}\vec J$ is
$L^2_{\text{loc}}(U)$ so $\tilde{\underline{v}}$ is well defined on
$Q$, as a function.

Next, we extend the definition of $\tilde{\underline{v}}$ to all
of $\Bbb R^n\times \Bbb R$ by $0$:
$$
\underline{v}\equiv\cases\tilde{\underline{v}} \text{ on } Q\\ 0
\text{ on } Q^c\endcases.\tag 2.3
$$
Next, since $\nabla_{(x, t)}\psi= (\nabla \psi, \partial_t \psi) = 0$
a.e. on $Q^c$ it follows that
$$\vec J = \rho \underline{v} \text{ a.e. on } \Bbb R^n \times \Bbb
R,\tag 2.4
$$
and moreover we can compute on $Q$, where $\frac{\psi}{|\psi|}$ is
bounded and continuous as a map from $Q$ to $S^1$:
$$
\aligned
\nabla \frac{\psi}{|\psi|} &= \eta^{-1}\nabla \psi - \eta^{-2} \psi
\nabla|\psi|\\
&=\eta^{-1}\nabla \psi - \eta^{-2} \psi  \frac{1}{2}  \eta^{-1}(\bar \psi
\nabla \psi + \psi\nabla\bar \psi)\\
&= \eta^{-1} \bar \psi^{-1} \frac{1}{2} [\bar\psi \nabla \psi - \psi
\nabla \bar \psi]\\
\nabla \psi &= \nabla (\eta\frac{\psi}{|\psi|}) = \bar
\psi^{-1}\frac{1}{2}[\bar \psi \nabla \psi - \psi \nabla \bar\psi]
\\
&+\frac{\psi}{|\psi|}\nabla\eta\\
\frac{\bar\psi}{|\psi|} \nabla \psi &= \eta^{-1} \frac{1}{2} [\bar\psi
\nabla \psi - \psi \nabla \bar \psi]+ \nabla \eta\\
&= i \eta \underline{v} + \nabla \eta.
\endaligned
$$
For each $t$,
$$
\int_{\Bbb R}|\nabla \psi|^2 = \int_{N^c_t}|\nabla \psi|^2 =
\int_{N^c_t} |\frac{\bar \psi}{\psi} \nabla \psi|^2 = \int_{N^c_t}
|\nabla \eta |^2 + \int_{N^c_t} \eta^2|\underline{v}|^2\tag 2.5
$$
$$
= \int_{\Bbb R} |\nabla \eta|^2 + \int_{\Bbb R}\eta^2|\underline{v}|^2
$$
since $\nabla \eta$ and $\eta\underline{v}$ are zero a.e. on $N_t$.
We therefore have that $\underline{v}\in L^2(\eta^2 dx)$.

Next, we construct the angular velocity function $\omega$, in a
similar way.

Again, on $Q$ we have
$$
\partial_t \frac{\psi}{|\psi|} = \frac{1}{|\psi|}
\dot\psi-\frac{\psi}{|\psi|^2} \partial_t |\psi|
$$
$$
= |\psi|^{-1}[\dot\psi - \frac{\psi}{|\psi|}
\frac{1}{2}|\psi|^{-1}(\bar\psi\dot\psi + \Dot{\Bar\psi}\psi)]
$$
$$
= |\psi|^{-1}[\frac{1}{2} \dot\psi - \frac{1}{2}
\frac{\psi}{|\psi|^2}{\Dot{\Bar\psi}}\psi]
$$
$$
=\eta^{-1}\frac{1}{2} [\dot\psi - \frac{\psi}{\bar\psi} \Dot{\Bar
\psi}]
$$
$$
\frac{1}{2} \eta^{-1} [\bar\psi\dot\psi-\psi\Dot{\Bar\psi}]\frac{1}{\bar\psi}.
$$
Let, on $Q$:

$$
\omega \equiv\frac{\bar\psi}{|\psi|}\partial_t \frac{\psi}{|\psi|} =
\eta^{-2} \frac{-i}{2} [\bar\psi(i\dot\psi) +
\psi(\Dot{\Bar{i\psi}})]\tag 2.6
$$
$$
= \eta^{-2} (-\frac{i}{2})[\bar \psi(-\Delta \psi) + \psi(-\Delta\bar
\psi)] + \eta^{-2} O(F(|\psi|)).
$$
Then, since $\eta^{-2}\bar\psi$ is continuous in $Q$, and
$$
-\Delta\psi \in H^{-1}(\Bbb R)
$$
$\omega$ is a distribution, for each $t$.

Also, since $\eta^2$ is bounded, continuous, $\nabla \eta \in L^2$,
$$
\eta^2 \omega \in H^{-1} (\Bbb R).
$$
Moreover, $\omega$ is trivially a distribution in $H^{-1}(Q)$.

\proclaim{Proposition} (Joint Continuity)
  Under the conditions of section I on $\psi,
\psi(x, t)$ is jointly continuous in $(x, t)$ in one dimension.
\endproclaim

The same holds in $\Bbb R^3/\{0\}\times \Bbb R$.

\demo{Proof}
$$
\psi (x, t) = e^{i\Delta t}\psi_0 - i\int^t_0 e^{i\Delta(t-s)}
F(|\psi(s)|) \psi(s) ds\tag A1
$$
$$
\nabla \psi(x, t) = e^{i\Delta t} \nabla \psi_0 - i\int^t_0
e^{i\Delta(t-s)}[G(\psi)\psi + F(|\psi|)\nabla \psi]ds
\tag A2
$$
$$
|\psi(x, t) - \psi(x,t')| \leq c \| \psi(t) - \psi(t')\|_{H^{1/2
 +\varep}}. \tag A3
$$

Therefore, to prove continuity, it is sufficient to prove continuity
in $t$ into $H^1(\Bbb R)$.  For this we use (A2).  Since $\nabla
\psi_0\in L^2$, the Spectral and Von-Neumann theorem gives continuity
in $t$ of the first term (in $L^2$).

As for the second term, continuity follows if the integrand in A2 is
$L^1(L^2(\Bbb R), dt)$, that is, if
$$
\|G(\psi)\psi + F(|\psi|)\nabla\psi\|_{L^2(\Bbb R)}\in L^1(dt).\tag 2.7
$$
Since $\psi$ is bounded and $\nabla \psi \in L^2$, the result
follows.  Similar computations work in the radial  case in three
dimensions.  \qed
\enddemo

In the $\eta, \underline{v}$ representation the IWC becomes:

When $\rho \neq 0$:
$$
i\hat{r}\cdot (\psi\nabla \bar\psi - \bar \psi \nabla \psi) =
2\rho
 \hat r\cdot \underline{v}  =  2\eta^2 \hat r \cdot \underline{v}.\tag 2.8
$$
When $\rho = 0,  \psi \nabla \bar\psi - \bar\psi \nabla \psi = 0$ and
hence we still have
$$
i\hat r \cdot(\psi\nabla\bar \psi - \bar\psi \nabla \psi) = 2
\eta^2
\hat r \cdot \underline{v} (= 0).\tag 2.9
$$
so that
$$\text{for all $x, t $\quad IWC: } \,
\eta^2 \hat r\cdot \underline{v} \leq 0.
\tag 2.10
$$
Furthermore, the $H^1$ bound gives
$$
\align
&\sup_t\|\nabla \eta\|^2_{L^2} \leq \Cal E_\rho < \infty\tag 2.11a \\
&\sup_t\|\eta|
\underline{v} |\|^2_{L^2} \leq \Cal E_{\underline{v}}< \infty
\tag 2.11b
\endalign
$$
using equation (8).

We can now use Sobolev estimates to conclude that

\proclaim{Proposition 2.1}  For $n=1, \eta $ is bounded and
continuous function of $x$.  For $n=3, \eta(r)$ is bounded and
continuous, for $r\neq 0$.  The bounds are uniform in time, $t$.

\endproclaim

The Schr\"odinger equation gives:
$$
\eta\dot\eta = \frac{1}{2} \partial_t (\eta^2) =
\frac{1}{2}\partial_t\rho = 2 \eta\eta'\cdot
\underline{v}
+ \eta^2\nabla \cdot \underline{v}.\tag 2.12
$$

Using this, we have
$$
\frac{1}{2} \partial_t\int_{R_1 \leq r \leq R_2} \eta^2 d^n x =
\int_{R_1\leq r \leq R_2} [\nabla (\eta^2)\cdot
\underline{v}+
\eta^2\nabla\cdot
\underline{v}
]d^n x = \int_S\eta^2 \underline{v}
\cdot d\overarrow S
$$
where $S$ is the closed surface of the domain $R_1\leq r \leq
R_2$. For $n=1$ , $ S$ is $\{ R_1, R_2\}$.  Let $\hat N$ be the
unit normal vector to surface $S$.  Using radial symmetry we
arrive at

\proclaim{Lemma 2.1}

a)  For almost all $t$,

$$
\align
\frac{1}{2} \partial_t\int_{R_1 \leq r \leq R_2}
\rho d^n x &= c_n(r^{n-1} \eta^2(r,t)
\underline{v}
\cdot \hat N
\bigg|^{R_2}_{r=R_1}\\
&= c_n(R_2^{n-1} \eta^2(R_2, t)
\underline{v}
\cdot \hat N(R_2, t) -
R^{n-1}_1\eta_2(R_1, t)
\underline{v}
\cdot\hat N(R_1, t)).
\endalign
$$

b) The above formula holds for all $t$, almost everywhere in $R_i, i =
1, 2$.
\endproclaim

\proclaim{Lemma 2.2 (Flux at infinity)}

Under our assumptions, the flux $\eta^2
\underline{v}
\cdot \hat N$
vanishes at infinity.
\endproclaim

\demo{Proof}  Since $\psi \in H^1$, we know that
$$
\int \eta^2 |
\underline{v}
|^2 d^n x < \infty \text{ which implies the
existence of a sequence $r_m \to \infty$, such that}
$$
$$
r_m^{n-1} \eta^2(r_m) |
\underline{v}
(r_m)|^2 \to 0 \text{ as } r_m \to
\infty,
$$
for each $t$ fixed.

Therefore, for each time $t$ fixed,

$$r_m^{n-1} \eta
^2(r_m) |
\underline{v}
(r_m) |\to 0, \tag 2.13
$$
otherwise $|\underline{v}
(r_m) |\to 0$.  But then, using the radial
symmetry and the fact that $\eta \in H^1_{\text{ radial}}$ it follows
that
$$
r^{n-1}\eta^2\leq c, \text{which also implies that}
$$
$$
r^{n-1}_m \eta^2(r_m) |
\underline{v}
(r_m) |\to 0.
$$\hfill\qed
\enddemo

\proclaim{Lemma 2.3 (Flux at zero)}  There is no flux through the origin
\endproclaim

\demo{Proof}   By the previous Lemma,
$$
\eta^2
\underline{v}
\cdot \hat N \big|^\infty_{r=R} = 0 -
\eta^2
\underline{v}
\cdot \hat N\big|_{r=R} \leq 0.\tag 2.14
$$

Hence,
$$\partial_t \int_{0\leq r \leq \infty} \eta^2 d^n x = 0 =
\eta^2
\underline{v}
\cdot \hat N\big|^\infty_{R\to 0} = -
\lim_{R\to 0} \eta^2
\underline{v} \cdot \hat N = 0.
$$
\hfill\qed
\enddemo
From now on, on a sphere of radius $R$, we denote
$$
\underline{v}\cdot \hat N(R,t) \equiv v_n (R, t).
$$

\proclaim{Lemma 2.4}  For each $R, \eta^2 \hat r \cdot v
 (R, t) \in L^1(dt)$,
uniformly in $R$.
\endproclaim

\demo{Proof}
$$
\frac{1}{2} \int^T_0 dt (\partial_t\int_{|x|\leq R} \eta^2d^n x) =
\frac{1}{2} \int_{|x|\leq R} \eta^2 d^n x \big|^T_{t=0} = \int^T_0
\eta^2
\hat r\cdot  \underline{v}
(R, t) dt \tag 2.15
$$
and since the last integrand is positive, by IWC, it follows that
$$
\int^T_0 |\eta^2
\hat r\cdot \underline{v}
(R, t)|dt = \frac{1}{2}\left(\| E ( |x| \leq
R)\psi(T)\|^2_2 -\|E(|x|\leq R) \psi(0) \|^2_2 \right)< \infty
$$
since the left hand side is monotonic increasing in $T$, the result
follows:
$$
0\leq \int^T_0\eta^2
\hat r\cdot \underline{v}
(R, t) dt \leq \frac{1}{2}
\|\psi_0\|^2_2.\tag 2.16
$$
In particular, $\eta^2 v_n (R, t) \to ''0'' $ as an $L^1$
function.\qed
\enddemo

\proclaim{Proposition 2.6}  $$\text{Weak}-\lim_{t\to\infty} \eta(t)
\equiv u, \text{exists in } H^1.\tag 2.17
$$
\endproclaim

\demo{Proof}
$$
\int_{R_1\leq |x|\leq R_2} \eta^2 d^n x \big|^T_{t=0} = 2 \int^T_0
dt \left(\eta^2 \hat r \cdot \underline{v} \big|^{R_2}_{r=R_1}
\right.\tag 2.18
$$
For each $R_1, R_2$, the r.h.s. converges as $T\to \infty$, by
Lemma (2.4).  Hence, for each interval $I\subset \Bbb R$
$$
\lim_{t\to\infty} \int \chi_I(r)\eta^2 d^n x
$$exists for $\chi_I$ the characteristic function of $I$.
\enddemo

The finite linear combinations of such $\chi_I, \sum \alpha_i
\chi_{I_i}$ are dense in $L^p, 1 \leq p < \infty$.

Hence, given $\varep > 0$, and $\varphi \in L^p, 1\leq p <
\infty$, we can find an element $\varphi_\varepsilon = \sum^N_{i=1}
\alpha^\varepsilon_i \chi^\varep_{I_i}, N<\infty$,
and $\|\varphi - \varphi_\varep \|_{L^p} \leq \varep$.

So, $$\lim_{t\to\infty} \int \varphi\eta^2 d^n x = \lim_{t\to\infty}
\int(\varphi-\varphi_\varep)\eta^2 d^n x
$$
$$
+\lim_{t\to\infty} \int \varphi_\varep \eta^2 d^n x.
$$

The first term on the r.h.s. is bounded by
$$
\|\varphi- \varphi_\varep \|_p \sup_t\| \eta^2\|_{L^{p'}}\leq C
\varep
$$
for all $p'\leq \frac{n}{n-2} ( p' \leq 3 \text{ for } n=3, \infty
\text{ in } 1-\dim) $ and the second term converges as $t\to
\infty$.

Hence, we have

\proclaim{Lemma}  If $\lim_{t\to \infty} \int\chi_I (r)\rho_t d^n x$
exists for all intervals $I$, then $\rho_t\overset{\omega}\to{\to}\rho_\infty$
in $L^p$, for any $p>1$, provided
$$
\sup_t \|\rho_t\|_{L^p} < C < \infty.
$$
In particular, in our case, since $\eta\in H^1$,

$$
\align
\eta^2\in L^q, & 1 \leq q \leq 3 \text{ for } n=3\\
               & 1 \leq q\leq \infty \text{ for } n=1,
\endalign
$$
it follows that
$$
\omega-\lim_{t\to\infty} \eta^2_t \overset{L^q}\to{\to} u^2 \geq 0,
u^2 \in L^q, q > 1.\tag 2.19
$$
\endproclaim

Now let $\varphi\in C^\infty_0$ and the dimension to be 1.  For all
$t, t'$, we have
$$
\int\varphi \nabla(\eta^2_t - \eta^2_{t'}) = -
\int(\nabla\phi)(\eta^2_t - \eta^2_{t'}) \to 0 \text{ as }
t,t'\to\infty
$$
since $\eta^2_t$ is Cauchy in $L^2$, say.

Since $C^\infty_0$ is dense in $L^2$ and $\sup_t\|\nabla
\eta^2_t\|_2 < \infty, \nabla\eta^2_t$ is Cauchy and hence
converges in $L^2$, weakly.  Here we used that $\|\nabla
\eta^2_t\|_2 \leq 2\|\eta_t\|_\infty \|\nabla \eta_t\|_2 \leq
C\|\nabla \eta_t\|^2_2. $ In three dimensions a similar argument
applies:
$$
\nabla (\eta^2_t) \in L^{3/2}, \text{ therefore }
$$
$$
\eta^2_t \to u^2\text{ weakly in } W^{1,3/2}.
$$
Next, using that
$$
\sup_t\int |x|^2\rho_t d^n x < \infty
$$
so that $\{\eta^2_t\}_{t\geq 0}$ is precompact in $H^s(\Bbb R) ,
0 \leq s < 1$ and $(n=3)$ in $W^{s, 3/2}(\Bbb R^3)\quad 0 \leq s<
1$.

Using that $W^{1, 3/2} \hookrightarrow H^s$ for $s = 1/2 (n=3)$.

We conclude
\proclaim{Proposition}
$$
\eta^2_t\to u^2\text{ strongly in } H^s(\Bbb R)\quad 0\leq s < 1
$$
and in $W^{s, 3/2}, 0\leq s <1, H^s(\Bbb R^3),  0\leq s  < 1/2$ in
3 dimensions.

Applying Sobolev embedding theorem again:
$$
\eta^2_t \to u^2\text{ strongly in } L^p(\Bbb R^n), p< \frac{n}{n-2}
$$
$n\geq 3,$ and all $p$ in $1-\dim$.
\endproclaim
Now, we use that to prove {\bf weak} convergence of $\eta_t$ in
$L^p$.

Let $\varphi \in C_0^\infty$, supported in some compact interval
$I\subset \Bbb R,$

$$
|I|< \infty.
$$
Then,
$$
J \equiv |\int \varphi(\eta_t - u) |\leq \int |\varphi|\chi_1
|(\eta_t - u) (\eta_t + u) | (\eta_t + u)^{-1}
$$
$$
+ \int|\varphi|\chi_2 | | \eta_t - u|
$$
where $\chi_1$ is the characteristic function of the set where $\eta_t
+ u\geq \varep > 0$, and $\chi_2\equiv 1-\chi_1$.

So,
$$
J\leq \frac{1}{\varep}\|\eta_t^2 - u^2\|_{L^p}
\|\varphi\|_{L^{p'}} + \varep |I|^{1/p} \|\varphi\|_{L^{p'}}
$$
where we used that
$$
\chi_2|\eta_t - u|\leq \chi_2(\eta_t + u) \leq \varep \chi_2
$$
since $\eta_t, u\geq 0$.

Given $\varep_0 > 0$, choose $\varep$ s.t.
$$
\varep|I|^{1/p} \|\varphi\|_{L^{p'}}< \varep_0/2.
$$
For this $\varep$, choose $T(\varep\varep_0)$ s.t. for all $t>
T(\varep\varep_0)$
$$
\|\eta^2_t - u^2\|_{L^p} \|\varphi\|_{L^{p'}}\leq \varep\varep_0/2.
$$
Hence $J<\varep_0$.  This can be done for all $2\leq p <q_*$
$$
\align
q_* &= 3      \text{ when } n=3\\
q_* &= \infty \text{ when } n=1.
\endalign
$$
Hence the above applies to all $1< p' \leq 2$ in $1-\dim$ and $3/2 <
p'\leq 2$ in $3-\dim$.

So, in $n=1, $ since $\eta_t$ is uniformly bounded in $L^p$ for all
$p\geq 2, \eta_t\overset{L^p}\to{\to} u$ weakly for $2\leq p <
\infty$.

In three dimensions, $\eta_t$ is uniformly bounded in $L^p$,
$$
2\leq p \leq 6
$$
so $\eta_t\overset{L^p}\to{\to} u$ weakly for $2\leq p <3.$
We now take, as before $J\equiv |\int \varphi\nabla
(\eta_t-\eta_{t'})|$ and argue as above, to conclude that
$$\nabla \eta_t\text{ is Cauchy in } L^2
$$
which implies that $\nabla \eta_t\overset{L^2}\to{\to}\nabla \tilde u$
weakly.  But then, for $\varphi\in
C^\infty_0 ,  \,  \int\varphi\nabla(\eta_t-\tilde u) = - \int\nabla
\varphi(\eta_t-\tilde u) \to 0$ implies that $\tilde u = u$, since
$\eta_t\overset{L^2}\to{\to} u$ weakly. \qed

As corollaries we have

\proclaim{Proposition 2.7}
$$
\eta_t\overset{s} \to{\to}u \text{ in } H^s, 0\leq s < 1.
$$
\endproclaim

\noindent{\bf PF}

Follows from weak convergence in $H^1$ and from the precompactness in
$H^s, s< 1.$

\proclaim{Proposition 2.9}
$$
u \in H^1 \, \quad\|u\|_{L^\infty(\Bbb R)}< \infty \quad u \in
\text{Lip} (\frac{1}{2}-\varep) \text{ in } \Bbb R.
$$
$ |u| \leq c/\sqrt{r} \text{ in } H^1_{\text{rad}}(\Bbb R^3)$.
 $u
\in \text{Lip} (\frac{1}{2})$ away from zero in $\Bbb R^3$.
\endproclaim

\head Section 3  Pointwise Convergence of $\eta,
\underline{v}$ and other
properties\endhead

We have seen that $u\geq 0$, is continuous in $\Bbb R$, and in $\Bbb
R^3/\{0\}$. Hence, if $u(x_0)> 0$, there is an interval $\tilde I$
around $x_0$ where $u$ vanishes for the first time at its end points,
in one dimension.  In three dimensions the same holds, provided the
origin is outside $r\in \tilde I_{x_0}$.  $r=|x|$.  If $I\subset
\tilde I_{x_0}$, $I$ away from the boundary, then
$$
u(x) > \delta > 0 \text{ for all  } x \in I.
$$

We will now analyze $u$ in such $I$.

\proclaim{Proposition 3.1}
$$
\eta_t \to u\text{ pointwise in }\Bbb R \text{ or }\Bbb R^3/\{0\}\text{
uniformly}.
$$
\endproclaim
\demo{Proof}
$$
|\eta_t (x) - u(x)|\leq C\|\eta_t - u\|_{H^{1/2+\varep}(\Bbb R)}\to 0
 \text{ as } t \to \infty
$$
by proposition  2.7.

Similar estimate holds in $\Bbb R^3_{\text{rad}}  / \{0\}.$

\enddemo

\proclaim{Lemma 3.2}   On $I$, the regularity of $\psi$ carries over
that of $\eta$.
\endproclaim

\demo{Proof}  Say $s = 1 + \mu, \mu < 1$.  Let $\psi = a+ i b$.
$$
\| D^s\rho \|= \| D^\mu(\bar \psi \nabla\psi + \psi\nabla \bar \psi) \|
\leq C\| D^\mu \psi\| \| \nabla \psi\| + C \| \psi\| \| D^{1+\mu} \psi \|
< \infty
$$
by assumption  (that $\psi\in H^s)$.

$\|D^{1+\mu}\chi_I|\psi | \|=\|D^\mu D(\chi_I\sqrt{a^2 + b^2}) \|$
$$
= \| D^\mu[(D\chi_I)|\psi|] + D^\mu\{ \chi_I(a^2 + b^2)^{-1/2} (a\nabla a
+ b\nabla b)\|
$$
so, we only need to bound
$$
\|D^\mu[(a^2+b^2)^{-1/2}\chi_I]\|_{L^q}
\text{ for some $q$ large enough}.
$$

This last expression is equal to
$$
\|D^{-1 + \mu} D[(a^2 + b^2)^{-1/2} \chi_I]\|_{L^q} \leq C
D[(a^2+ b^2)^{-1/2}\chi_I]\|_{L^{\tilde q}}
$$
for some $2 \leq \tilde q \leq 3$.
But $D[(a^2+b^2)^{-1/2} \chi_I] = - (a^2+ b^2)^{- 3/2}(a\nabla a + b
\nabla b) \chi_I + 0(1) \in L^{\tilde q}$ since $a, \chi_I (a^2 +
b^2)^{-3/2} \in L^\infty, \nabla a, \nabla b \in L^{\tilde q}$.\qed
\enddemo

\proclaim{Proposition 3.3}

Recall that $\rho \underline{v}  \to ''0''$ as $L^1(dt)$ for $t\to \infty$.

a)  On $I\subset \tilde I$ we then have: $(\rho\equiv \eta^2)$

$$
\eta\hat r\cdot
\underline{v}  \to ''0'' \text{ and } \hat r \cdot
\underline{v} \to ''0'' \text{ as }
L^1(dt).
$$

b) $\chi_{\tilde I} \eta\hat r\cdot
\underline{v} \overset{ L^2}\to{\to} 0$  strongly and
$\chi_I \hat r \cdot \underline{v}
\overset{ L^2}\to{\to} 0 $ strongly, for $t_n \to \infty.$
\endproclaim

\demo{Proof}

Part (a) follows, since on $I, \eta > \delta > 0$, uniformly in time,
since on $I, u > \delta' \geq 2 \delta$, and $\eta_t \to u$ pointwise
by Proposition 3.1.  To Prove (b), we note that $\|\eta\hat r \cdot \underline{v}\|_{L^2}
< \infty$ uniformly in $t$, since
$$
\eta^2 \hat r \cdot \underline{v} = \frac{i}{2} (\bar \psi \nabla \psi
- \psi \nabla \bar \psi) \cdot \hat r\text{ then }
\eta\hat r \cdot v = \vec J\cdot \hat r/\eta (\text{ on }Q)
$$
and
$$|\vec J/\eta|\leq |\nabla \psi| \in L^2,
$$
and $0$ otherwise, with uniform (in $t$) $L^2$-norm.
\enddemo

We now use the following Lemma.

\proclaim{Lemma}  If $\| f_n\|_{L^p} < C < \infty, f_n \to f$
pointwise, then $\chi_J(f_n-f) \to 0$ in $L^q$, for all $q< p$ and $J$
compact.
\endproclaim

By this Lemma the result (b)follows if we can prove that for a
sequence $t_n \to \infty$,
$$
\sup_n \|\chi_{\tilde I}\eta_{t_n}\hat r\cdot \underline{v} (t_n)
\|_{L^q} < \infty
$$
for some $q> 2$.

This follows from
Strichartz estimate on compact time intervals, as it implies that
$\|\nabla\psi_n \|_{L^3}\leq c\|\psi(0) \|_{H^1}$ uniformly in $n$,
and
$$
|\nabla \psi_n |^3 \geq |\eta_n \hat r \cdot \underline{v}_n |^3.
$$
\qed
\proclaim{Proposition 3.4}
$$
\chi_I\eta_n \to \chi_I u \text{ strongly in } H^s, s < 3/2.
$$
\endproclaim

\demo{Proof}The result follows from weak convergence and compactness
of $\chi_I \eta_n$ in $H^s, s < 3/2 $, if we can prove uniform
boundedness in $H^s,  s< 3/2.$

The weak convergence follows from $(\varphi \in C^\infty_0)$
$$
|\int \varphi \chi_I D^s(\eta_n - \eta_m)|=|- \int(\varphi\chi_I)'
 D^{s-1} (\eta_n - \eta_m)|
$$
$\leq C_\varphi \| D^{s-1} (\eta_n - \eta_m) \|\to 0$ as $ n,m\to
\infty$ by proposition 2.7 provided $0\leq s - 1< 1$.
So, the weak convergence will follow if $\|\chi_I \eta_n\|_{H^s} < c <
\infty$ uniformly in $n$.  For $s\leq 3/2$ this follows from Local
smoothing estimates on compact time intervals. \qed
\enddemo

\proclaim{Corollary 3.1}
$$
\chi_I u \in H^s, \quad s< 3/2
$$
$$
\chi_I u \in W^{1, 3}.
$$
\endproclaim

\head Section 4 - The equation for $u$\endhead

Next, we would like to rewrite the NLS in an equivalent form, in terms of $\eta,
\theta:$

Let us compute, formally:
$$
\align
\psi &= \eta e^{-i\theta}\quad \nabla \psi = (\nabla \eta) e^{-i\theta}
-i\theta'\eta e^{-i\theta} \tag 4.1 \\
i\dot{\psi} &= i\dot{\eta} e^{-i\theta} + \dot{\theta} \eta
e^{-i\theta}\quad F(|\psi|) = F(\eta)\\
-\Delta \psi&= - \nabla \cdot (
\eta' e^{-i\theta} - i\eta\theta' e^{-i\theta}) =- (\Delta \eta)
e^{-i\theta}\\
&+ 2 i \eta' \cdot \theta' e^{-i\theta} + \eta\theta^{\prime 2}
e^{-i\theta} + i\eta \theta^{''} e^{-i\theta}\endalign
$$
where $'$ stands for gradient, $''$ for $\Delta: \eta'\cdot
\theta' =\nabla \eta\cdot \nabla \theta, \theta^{''} = \Delta \theta$.

We have:

$$
\align
\dot{\eta}&=2\eta'\theta' + \eta\theta^{''}= 2\nabla\eta \cdot \nabla
\theta + \eta \Delta \theta  \tag 4.2a \\
-\Delta\eta &+ F(\eta)\eta = \dot\theta \eta - \eta \theta^{'2}.\tag4.2b
\endalign
$$

Next, we embark on a way to make sense of these equations [see also
AC, Ger, Gre, LZ] in some parts of the space time.
Let $I\subset \tilde I$ be such an interval where $u\geq \delta > 0$
for all $r\in I$.

Since $\eta \to u$ uniformly on $I$ (Proposition 3.4)  we have
that
$$
\eta_t \geq \delta / 2 > 0,\quad\forall\,   t > T_0, \text{ on } I.\tag 4.3
$$
We therefore have that there exists a (good) box $B$ in space-time,
where $\eta\geq \delta / 2 > 0$ for all $(x, t)\in B$,
$$
B=\{ x\in I,\quad t\geq T_0\}.
$$

$I$ is any interval strictly between two consecutive zeros of $u$.
We always assume $I$ does not contain the origin in $\Bbb R^3$.
We also recall the notation
$$
\aligned
Q&= \{(x,t) \in \Bbb R^n \times \Bbb R, \eta(x, t)\neq 0\}\\
N_t &= \{ x\in \Bbb R^n, \eta (x, t) = 0\}.
\endaligned
$$
Since $\eta$ is continuous in $\Bbb R$ (or $\Bbb R^3/ \{ 0 \}$), the
set $N_t$ is closed, $N_t^c$ is open.

Hence the function
$$
\frac{\psi(x, t)}{|\psi(x, t)|} = \frac{\psi(x, t)}{\eta(x, t)} : Q\to
S^1\tag 4.4
$$
is continuous (jointly in $(x, t)$).

For each $t, N_t^c$ is a collection of open connected disjoint sets,
$N_{t,i}^c  i = 1, 2, \dots$ .

The key to constructing the phase function $\theta$, with the right
properties, is the existence of continuous extension to the (universal)
covering space of $S^1, \Bbb R$.

The conditions for these are known:

\proclaim{ Theorem } (Lifting Lemma)

Let $p: E\to B$ be a covering map of a topological space $B$, by $E$.

let $p(e_0)= b_0$.  Let $f:Y\to B$ be a continuous map of a
topological space $Y$ into $B$, with $f(y_0)= b_0$.  Suppose that $Y$
is path connected and locally path connected.  The map $f$ can then be
lifted to a map $\tilde f: Y\to E$ such that
$$
\tilde f(y_0) = e_0, \, \text{ if and only if}\tag 4.5
$$
$$
f_*(\pi_1(Y, y_0)) \subset p_*(\pi_1(E, e_0))\tag 4.6  (Condition
L)
$$
Furthermore, if such a lifting exists, it is unique.

Here $\pi_1(A, a)$ denotes the first fundamental group of loops of $A$
starting at the point $a$.

$f_*, p_*$ are the induced maps on the loop spaces by $f$ and $p$,
respectively.
\endproclaim
Using this theorem [GH] we can conclude the following :

$\bullet$  Since in
general the space $Q$ is not simply connected, and since $\pi_1(\Bbb
R, e_0)$ is trivial, condition (L) will not be satisfied in general.
So, the lifting of
$$f\equiv \frac{\psi}{\eta}: Q\to S^1
$$
does not exist.

$\bullet$  For each $t,$ and a connected part of $N^c_t$, a lifting $\Cal Q_i(x,
t)$ exists and is unique; but it is {\bf not} continuous in $t$.

$\bullet$  A unique lifting $\theta_B(x, t)$ which is jointly continuous exists on
  good boxes $B$, since boxes have trivial fundamental group.  It is
  this $\theta_B(x, t)$ which we will use from now on.

Using $\theta_B\equiv \theta$ we have on $B$:
$$
\frac{\psi}{|\psi|} = e^{-i\theta (x, t) }, \,  (x, t) \in B\tag 4.7
$$
with $\theta (x, t)$ jointly continuous in $(x, t)$.

Uniqness is guaranteed, up to a choice of reference point, namely we
choose a point $(x_0, T_0)$ where $\frac{\psi}{|\psi|}$ is not on
$\Bbb R^-$, and we choose $\frac{\psi}{|\psi|}$ to be on the first
Riemann sheet of the $\ell n z$ function (with cut on $\Bbb R^-$).

We therefore established the following crucial identities
$$
\psi = \eta
 e^{-i\theta (x, t) }, \, (x, t)\in B\tag 4.8a
$$
$$
\ell n \frac{\psi}{\|\psi\|} = - i\theta.\tag 4.8b
$$

\remark{Remark}  This construction is limited to the case of $\psi$
being continuous.
\endremark

\proclaim{Corollary}  Derivatives commute with lifting. Hence
$$
(f')_* = (f_*)'.\tag 4.9
$$

\endproclaim

To take a derivative of $\psi$ requires the use of the chain rule
(which uses the above Corollary).
$$
\frac{\partial \psi}{\partial x} = (\partial_x\eta) e^{-i\theta}+-
i\eta(\partial_x\theta) e^{-i\theta} \tag 4.10
$$
and similarly for $\frac{\partial\psi}{\partial t}$.

Once the above formula is established,
it follows from the fact that $\nabla_x \psi, \partial_x \eta$ are
in $L^2$, that $\nabla_x\theta$ is also in $L^2$ since $\eta \geq
\delta /2$ in $B$.

For $\dot\theta$, we use that $\dot\theta$ is well defined in the
sense of distributions, hence
$$
\int_B\int g(t) \varphi(x) \dot\theta d x dt = - \int_B\int
g'(t)\varphi(x) \theta (x, t) dxdt.
$$

Next, we compute higher derivatives, in the following weak sense: test
functions $\varphi \in H^1_0(I)$.

Since $\frac{\partial\psi}{\partial x} \in L^2$, the above terms can
be differentiated in the sense of distributions.

$$
\int\varphi\frac{\partial}{\partial x}\left[ \eta' e^{-i\theta} -
i\theta'\eta e^{-i\theta} \right] = \int \varphi \Delta \psi =
$$
$$
= \int\varphi\left[ \eta^{''} - i \theta^{''}\eta- i\theta'\eta' -
i\theta'\eta' - \theta^{'2}\eta\right] e^{-i\theta}.\tag 4.11
$$
In fact, each term makes sense separately since $e^{-i\theta}$ is
bounded, cont, $\eta^{''}, \theta^{''}$ are derivatives of locally
(square) integrable functions and so are distributions;
$$
\theta^{'2}\eta_{|_I}=\eta^{-1}\theta^{'2}\eta^2_{|_I} \leq
(\delta/2)^{-1}(\eta\theta')^2\in L^1_{\text{ loc}}(I).
$$
Similarly, $\theta'\eta'|_I \leq (\delta/ 2)^{-1}(\eta\theta')\eta' \in
L^1_{\text{loc}} (I)$ since both $\eta'$ and $\eta\theta'$ are in
$L^2_{\text{loc}}$.

In particular
$$
\bigg|\int\varphi \theta^{''}\eta\bigg| = \bigg| - \int\varphi'\theta'\eta - \int\varphi\theta'\eta'\bigg|
$$
$$
\leq \|\varphi\|_{H^1(I)}\left[\|\eta\theta'\|_{L^2(I)} +
\left(\frac{\delta}{2}\right)^{-1} \|\eta'\|_{L^2}\|\eta\theta'
\|_{L^2(I)}\right].
$$

On $B$, we have
$$
\int\varphi i\frac{\partial \psi}{\partial t} = \int i
\varphi\frac{\partial \psi}{\partial t} = \int i \varphi \dot \eta
e^{-i\theta} - 1\int \varphi\eta\dot\theta e^{-i\theta}.\tag 4.12
$$

Let $\varphi \in C^\infty_0(\tilde I)$.  Then

\proclaim{Lemma 4.1}
$$
(\varphi, -\Delta \eta) = ( -\Delta \varphi, \eta) \to (-\Delta
\varphi, u)\tag 4.13
$$
as $t\to \infty$.

This follows since $\eta\overset{ L^2}\to{\to} u$.

Similarly
$$
(\varphi, F(\eta) \eta) \to (\varphi, F(u)u).\tag 4.14
$$
\endproclaim

\proclaim{Lemma 4.2}
$$
\frac{1}{T}\int^T_{T_0} (\varphi, \eta |\nabla \theta|^2) dt \leq
c_\delta |I|^{1/2} \|\psi\|_2 \|\psi\|_{H^1} T^{-1/2} .\tag 4.15
$$
\endproclaim

\demo{Proof}
$$
\int^T_{T_0}(\varphi, \eta \theta^{'2}) dt = \int^T_{T_0}\int_I
\varphi|\nabla \theta|^{3/2} \eta|\nabla \theta|^{1/2} d x dt
$$
$$
\leq \left(\int_I dx \int^T_{T_0} |\varphi|^2 |\nabla \theta |^3
dt\right)^{1/2} \left(\int_I d x \int^T_{T_0} \eta^2 |\theta'|
dt\right)^{1/2}\tag{*}
$$
The second double integral is bounded as $T\to \infty$ by $(I
\sup_{r\in I} \int^\infty_0 |\eta^2\theta'| dt)^{1/2} \leq
|I|^{1/2} \|\psi\|_2$ by Lemma 2.6

Note that the proposition only gives a bound on the radial derivative
$\hat N \cdot \nabla \theta$, but for radial functions it is all we
need here.

The first double integral in (*) is bounded by the square root of
$$
\frac{8}{\delta^3} \int^T_{T_0} dt \int_I |\varphi |^2 \eta^3|\nabla
\theta|^3 dx \leq \frac{8}{\delta^3} \int^T_{T_0} \int \chi |\nabla
\psi|^3dx dt
$$
$$
\leq \frac{1}{\delta^3} T^{1/4}  \left(\int^T_{T_0} \|\nabla
\psi\|^4_3 dt\right)^{3/4}.
$$
By Strichartz inequality
$$
\int^{T'+1}_{T'}\|\nabla \psi\|^4_3 dt \leq c\|\psi\|^2_{H^1}, \tag
4.16
$$
so, the r.h.s. is bounded by
$$
\lesssim T^{1/4}  T^{3/4} = T .
$$
Hence, the first double integral in (*) is bounded by $T^{1/2} $;
so
$$
\int^T_{T_0} (\varphi, \eta\theta^{'2}) dt \leq c_\delta |I|^{1/2}
\|\psi\|_{H^1} \|\psi\|_{L^2} T^{1/2}.
$$\qed
\enddemo

Formally, we have:
$$
\int_{\Bbb R^n} \eta^2\dot \theta d^n x = \int|\nabla \eta|^2 + \int
F(\eta^2)\eta^2 + \int\eta^2|\nabla \theta|^2 \leq C < \infty\tag 4.17
$$
$$
\frac{1}{T} \int_0\int_{\Bbb R^n}\dot\theta \eta^2 d^n x dt =
\int_{\Bbb R^n} \frac{\theta(r, T)}{T} \eta^2(r, T) d^N x -
\frac{1}{T} \int_{\Bbb R^n} \theta(r, 0) \eta^2(r, 0) d^n x
$$
$$
-\frac{2}{T}\int^T_0 \int_{\Bbb R^n} \theta \eta(2\nabla \eta \cdot
 \nabla \theta + \eta \Delta \theta) d^n x dt.\tag 4.18a
$$
The last term is equal to
$$
\frac{2}{T} \int^T_0 \int_{\Bbb R^n} \eta^2|\nabla \theta|^2 d^n x dt
- \frac{2 C_n}{T} \int^T_0 \left(r^{n-1} \eta^2 \theta\theta'\big|
^\infty_{r=0}\right)dt.\tag 4.18b
$$

\demo{Proof}  To get (4.17) we multiply equation (4.2b) by $\eta$
and integrate over all space.

The last inequality follows from the energy estimate.  To get
(4.18a) we average over time multiplied by $\eta$, and integrate
by parts in $t$: that, together with (12) gives (18a).

To get (4.18b) we integrate by parts again:

$$
\int^T_0 \int_{\Bbb R^n}\theta \eta(2 \nabla \eta\cdot \nabla \theta +
\eta \Delta \theta) d^n x dt = \int^T_0\int_{\Bbb R^n} \theta (\nabla
\eta^2\cdot \nabla \theta + \eta^2 \nabla \cdot \nabla \theta) d^n x
dt
$$
$$
=\int^T_0\int_{\Bbb R^n} \theta \nabla (\eta^2\nabla \theta) d^n x dt
=- \int^T_0 \eta^2 |\nabla \theta |^2 d^n x dt + \int^T_0 (C_n r^{n-1}
\eta^2\theta\theta'\big|^\infty_{r=0} dt.$$\qed
\enddemo

Next, we use the above heuristics, to estimate the time average of
$\dot\theta\eta$.  We repeat the above argument with $\varphi \eta$
as a weight;

Then:
$$
\frac{1}{T}\int^T_{T_0} \int\varphi \dot \theta\eta^2 dxdt =
\int\frac{\theta(r, T)}{T}\varphi \eta^2_T dx -
\frac{1}{T}\int\theta(r, T_0) \varphi \eta^2_0 dx
$$
$$
-\frac{2}{T} \int^T_{T_0} \int \varphi \theta \eta (2\nabla\eta \cdot
 \nabla \theta + \eta\Delta\theta) dx dt \tag 4.19
$$
which is obtained by integration by parts in the $t$-variable, Recall
the notation $\eta_T = \eta(r, T)$.

Now, using that
$$
\eta(2\nabla\eta \cdot \nabla \theta + \eta \Delta\theta) = \nabla
\cdot(\eta^2\nabla \theta)
$$
and integrating by parts in $x$, the last term on the r.h.s of
equation (4.15) is equal to
$$
\frac{2}{T} \int^T_{T_0} \int(\nabla\varphi) \cdot\theta\eta^2 \nabla
\theta d x dt + \frac{2}{T}\int^T_{T_0} \int \varphi \eta^2
|\nabla\theta|^2 dx dt
$$
$$
-\frac{2}{T} \int^T_{T_0} (\varphi\theta \eta^2\nabla \theta \cdot \hat N
 \big|^\infty_{r=0} dt. \tag 4.20
$$
So, the l.h.s. of (4.18) is equal to a sum of 5 terms.

Since $\varphi$ is compactly supported the last, boundary term
vanishes.

For $t=T_0, \theta(r, T_0)$ is bounded for $r$ bounded, so the second term
on the r.h.s. vanishes like $O(\frac{1}{T})$ as $T\to \infty$.

By Lemma (4.2 ) the second term of equation (4.20) vanishes like
$T^{-1/2}$ as $T\to \infty$.

We are left with two terms: $\frac{\theta(r, T)}{T}$ term and
$\varphi'$ term.

We deal first with the $\frac{\theta(r,T)}{T}$ term:

\proclaim{Proposition 4.3}  For $r\in I \subset \tilde I,$
$$
\lim_{t\to\infty} t^{-1}\theta(r,t) \equiv - \Cal E,\text{  } \Cal
E \geq 0.\tag 4.21
$$
We have
$$
t^{-1}\theta(r,t) =t^{-1} \tilde{\Cal E}(t) + 0 (t^{-3/4}
 )\quad   t >> 1 \tag 4.22
$$
for some function  $\tilde{\Cal E}(t)$ {\bf independent} of $r$.
\endproclaim

\demo{Proof}
$$
t^{-1}\theta(r,t) = t^{-1} \theta(r_0, t) + t^{-1} \int^r_{r_0} \hat
r\cdot \nabla \theta(s, t) ds
$$
$$
\equiv t^{-1} \theta(r_0, t) + t^{-1} R (t)\tag 4.23
$$
and we choose:  $r_0 \in I$.
\enddemo

Now, for $r,r_0 \in I$
$$
|R(t)| \leq |r-r_0|^{1/2}\left(\int_I |\theta'(s, t)|^2
 ds\right)^{1/2} \leq |I|^{1/4} \delta^{-1} O(t^{1/4} )
$$
by Lemma 4.2

Hence (4.23) can be written as
$$
t^{-1}\theta(r,t) = t^{-1} \tilde{\Cal E}(t) + ( t^{-3/4} )\quad
t >> 1.
$$
with $\tilde{\Cal E}(t) \equiv \theta(r_0, t)$.

This proves (4.22).

We know that for $\varphi \in H^1(I),$
$$
\int \dot\theta\eta\varphi= \int\varphi(-\Delta \eta) + \int \varphi
F(\eta) \eta+\int \varphi \eta \theta^{'2}.\tag 4.24
$$

For $\varphi\in C^\infty_0(I)$, as $t\to \infty$ the r.h.s. converges
to
$$
\int (-\Delta\varphi) u + \int\varphi F(u) u + 0
$$
in the time mean.

The first two terms converge pointwise by Lemma (4.1), and the last
term converge in the mean (and therefore also for sequences) by Lemma
(4.2).

Let us apply $\frac{1}{T}\int^T_S(\varphi, \cdot) dt$ to equation
(4.2b):
$$
\frac{1}{T}\int^T_S(\varphi, -\Delta \eta) dt+
\frac{1}{T}\int^T_S(\varphi, F(\eta) \eta) dt + \frac{1}{T} \int^T_S
(\varphi, \eta\theta^{'2}) dt = \frac{1}{T} \int^T_S (\varphi, \dot
\theta\eta) dt.
$$

We now use that if $g(t)\to g$ as $t\to \infty$, then
$$
\frac{1}{T}\int^T_0 g(t) dt = \frac{1}{T}\int^S_0 g(t) dt +
\frac{1}{T}\int^T_S g(t) dt = o \left(\frac{S}{T}\right) + g + \frac{\varep
|T-S|}{T}.
$$
As $S\to \infty$, we can take $\varep\to 0$, so that
$$
\frac{1}{T}\int^T_0 g(t) dt\to g.
$$
So, for $\forall S\geq T_0$,
$$
\lim_{T\to \infty} \frac{1}{T} \int^T_S(\varphi, \dot \theta\eta) dt
= (\varphi, F(u) u) + (-\Delta \varphi, u) \tag 4.25
$$
and using (4.15), (4.19) with $T_0\to S$ at the lower integral
limit (in $t$):
$$
\frac{1}{T}\int^T_S \int\varphi \dot\theta \eta^2 dx dt = \frac{\tilde{\Cal
E}(T)}{T}\int\varphi\eta^2(x, T) dx  + \frac{2\tilde{\Cal
E}(T)}{T} \int^T_S \int_I (\nabla\varphi)\cdot \nabla\theta \eta^2 dx
dt + \tilde R(t)\tag 4.26
$$
and we want:
$$ \tilde R (T) = O (T^{-\varep}) .
$$
$$
\align
\tilde R(T) &= - \frac{1}{T}\int\theta(r, S)\varphi\eta^2_S dx +
 \Cal O (T^{-1/2} ) \tag 4.27\\
&+ \int \varphi T^{-1} R(T)\eta^2(x,T) dx +
2\frac{1}{T}\int^T_S \int (\nabla \varphi)\cdot \nabla\theta
R(t)\eta^2(x,t) dxdt\endalign
$$
(4.27) is obtained as follows: the first term comes from the second
term of (4.26)
(with $0\to S$).

The second term comes from the second term of (4.19) and the
estimate of Lemma (4.2 ).

The third term (and the forth) are obtained by replacing $\theta(r,
t)$ with (4.23), and recall the definition $\theta(r_0, t)\equiv
\tilde{\Cal E}(t)$.  Using (4.21) or (4.23), we see that the first
term of $\tilde R(T)$ is bounded by
$$CT^{-1} T^{1/4}  \|\eta\|^2_{L^2} \to 0\text{ as }
T\to \infty, \text{ like } T^{-3/4} .
$$
Similarly, the third term is vanishing like
$$
T^{-3/4} .
$$
Now, if we can prove that
$$
\int^\infty_0\int |\nabla \varphi \cdot \nabla \theta\eta^2| dxdt <
\infty
$$
it will follow that the last term of (4.27) vanishes like
$T^{-3/4} $.

Furthermore, by taking $S$ large as we want, it will follow that the
second term on the r.h.s. of (4.26) vanishes as $S\to \infty$, {\bf
provided} we know $a$-priori that
$$
T^{-1}\tilde{\Cal E}(T) \text{ is uniformly bounded in }T.
$$

We extend (4.25) to $\varphi\in H^1(I)$.  To this end let $f\in
H^1(I)$ and $\varep > 0$ be given.

Then there exists $\varphi_{\varep_n}\in C^\infty_0$ such that
$$
\| f- \varphi_{\varep_n} \|_{H^1} \leq \varep/n, n = 1,2,3..
$$
$$
\big|\frac{1}{T}\int^T_S (\varphi_{\varep_n} - \varphi_{\varep_m},
 \dot\theta\eta) dt\big|\leq |( \varphi_{\varep_n} -
 \varphi_{\varep_m}, F(u) u ) |+ |(\nabla(\varphi_{\varep_n} -
 \varphi_{\varep_m}) , \nabla u)| + \varep(T)\tag 4.28
$$
and $\varep(T) \to 0\text{ as } T\to \infty$.

Hence
$$
\Bigg|\frac{1}{T}\int^T_S(\varphi_{\varep_n} - \varphi_{\varep_m},
\dot\theta\eta) dt    \Bigg| \leq \varep(T) + \|\nabla u \|_2 \|
\nabla (\varphi_{\varep_n} - \varphi_{\varep_m})\|_2+
$$
$$+
\min \{\|\varphi_{\varep_n} -
\varphi_{\varep_m} \|_6 \,  \| F(u)
u\|_{6/5}; \|F(u)u)\|_2 \| \varphi_{\varep_n} - \varphi_{\varep_m}
\|_2\}\quad (\text{ in } n=3).\tag {*}
$$

And in dimension 1, we can replace 6 by $\infty$ and $6/5$ by 1.
Hence the l.h.s. of (*) is bounded by
$$
\varep(T) + C\varep/n.
$$
If we now let $T\to \infty$, it follows that
$$
a_n \equiv \lim_{T\to\infty} \frac{1}{T} \int^T_0
(\varphi_{\varep_n}, \dot \theta\eta) dt
$$
is Cauchy, and hence has a limit $a$.

Moreover, by (4.28)
$$
a_n\to (f, F(u)u) + (\nabla f, \nabla u).  \tag 4.29
$$
So, the linear operator $L$
$$
\lim_{T\to \infty} \frac{1}{T}\int^T_0 (f, \dot\theta\eta) dt: f\in
H^1(I)\overset{L}\to{\to} \Bbb C
$$
defined on $C^\infty_0(I)$ and extended by the above limit to all
$f\in H^1(I)$ satisfies the bound
$$
|Lf|\leq |(f,F(u)u)|+|(\nabla f, \nabla u)|
\leq c\|f\|_{H^1} \|u\|_{H^1} \leq M \|f\|_{H^1}.
$$
So, it is continuous.  Hence the extension of (4.25) to $\varphi\in
H^1 (I)$ follows.

On $I, \varphi\in C_0^\infty(I)$ implies that $\varphi u \in
H^1(I),$ the extension of (4.25) to $H^1(I)$ implies that
$$
\frac{1}{T} \int^T_S \varphi \dot\theta \eta^2 dxdt \, \text{  has a finite
limit}.
$$

So, the l.h.s. of (4.26) has a finite limit as $T\to \infty,
\tilde R(T)\to 0$, and the sum of the two remaining terms on the
r.h.s. of (4.26) is larger than
$$
\frac{\tilde \Cal E(T)}{T} \left[\int \varphi u^2 dx -
\varep_T(S)\right]
$$
with $\varep_T(S) \to 0\text{ as } T, S \to \infty$.

So, $\frac{\tilde \Cal E(T)}{T}$ can not diverge as $T\to \infty$, and
hence it has a limit $-\Cal E$.

Hence, as $T\to \infty$
$$
\frac{1}{T}\int^T_{T_0}\int \varphi \dot\theta \eta^2 dxdt \to - \Cal E
\int \varphi u^2 dx \tag 4.30
$$

If we choose $\varphi = \varphi_0/u\quad \varphi_0 \in C_0^\infty (I)$
we have: $\nabla \varphi = u^{-1} \nabla \varphi_0 - \varphi_0 u^{-2}
\nabla u \in L^2$  since $u\geq \delta \text{ on } I$.

So, using such $\varphi$ we have finally that
$$
\lim_{T\to \infty} \frac{1}{T} \int^T_0 \int \varphi_0 \dot\theta \eta
dxdt = - \Cal E \int \varphi_0 u dx \tag 4.31
$$
$$
= (-\varphi_0, \Delta u) + ( \varphi_0, F (u^2)u)
$$
or in the weak sense on test functions $\varphi_0\in C_0^\infty(I)$:
$$
-\Delta u + F(u^2)u =- \Cal E u \quad r\in I.\tag 4.32
$$

\proclaim{ Theorem 4.4}  The limiting state is a soliton for all $x\in
\Bbb R^n, n=1,3$
\endproclaim

\demo{Proof}

A-priori, we only know that $u$ is a soliton between two consecutive
zeros of $\rho$.  But, due to elliptic regularity, one can not have
a soliton made of ``patches'' of excited solitons (and or zero) with
the same energy.
\enddemo

\head Section 5  More general Solutions \endhead

Solutions which are completely incoming are special.  In many cases,
they consist of the soliton itself!  A close look at the proof shows
however, that the key assumption of IWC is used to prove the
integrability of the (local) flux.

We can therefore obtain, by similar arguments the following more
general results.

\proclaim{Proposition 5.1}  Let $\psi_0$ be as before, radial, $n = 3$.
Assume the solution of the NLS has integrable incoming (and hence also
outgoing) flux on any surface around the origin:
$$
\int^\infty_0 dt \int_{S_R}\hat n \cdot (\bar \psi \nabla \psi -
\psi\nabla \bar \psi)_\pm dS < \infty.
$$
\endproclaim

Here $(\cdot)_\pm$ stands for the positive / respectively negative
part of
$$
(\bar\psi \nabla \psi - \psi \nabla \bar \psi)\cdot \hat n
$$
where $\hat n$ is the unit normal vector to the sphere $S_R$ of radius
$R$ around the origin.  $dS$ is the surface element on $S_R$.

Then

a) For any compact interval $I$,
$$
\| \, | \psi (t) | - u \|_{ H^s(I)}
\to 0 \text{ as } t\to \infty
$$
$0<s<1$, and $u$ is a soliton solution of the NLS.

b)  The phase of the solution, $\theta (x, t)$ locally converges to the
soliton energy, in the sense that
$$
\frac{\theta (r, t)}{t} \to E, \text{ for $r$-fixed.}
$$

Proving that the incoming/outgoing flux is integrable in time is not
easy in general.

However, it should be noted that it follows from the following weaker
assumption, more suitable for applications:

\proclaim{ Proposition 5.2}

Assume $\psi_0$ is radial, $n=3$.

Suppose that the solution of the NLS with this initial condition has
a monotonic decreasing incoming flux, up to $L^1(dt)$ convergent
part, on any sphere around the origin.  Then the incoming/outgoing
fluxes are absolutely integrable in time on any sphere.
\endproclaim

The proof of proposition 5.2 follows since the total flux is
integrable (not absolutely!) in time on any sphere due to $L^2$
boundedness and conservation.

Hence, by spherical symmetry, on any sphere, at any time the flux is
either incoming or outgoing.  If it is always incoming, by the
integrability of the total flux, the (absolute) integrability of the
incoming wave is immediate.  If the flux turns from incoming to
outgoing after a finite time, it stays outgoing for later times,
($+L^1(dt)$ terms) by the assumption of monotone decrease, up to
$L^1(dt)$.  Again, since the total flux is integrable, the purely
outgoing part $(+L^1(dt))$ is absolutely integrable.

\medskip

\noindent{\bf Acknowledgements}

I wish to thank I. Rodnianski for very important discussions.

This work is partially supported by NSF grant DMS-0501043

\enddocument
\end

enddemo\enddocument


\Refs
\widestnumber\key{CEL -Sog}

\ref
\key AC
\by  T. Alazard, R. Carles
\paper  WKB analysis for the Gross-Pitaevskii equation...'',
\jour preprint arXiv 0710.0816V1 and cited references
\vol
\yr
\pages
\endref
\medskip

\ref
\key GH
\by  M. Greenberg and J. Harper
\paper  Algebraic Topology, a first course, Revised
\jour
\vol
\yr 1981
\pages
\endref
\medskip



\ref
\key Ger
\by P. G\'erard
\paper Remarques sur l'analyse semi-classique de l'\'equation de
Schr\"odinger non lin\'eaire
\jour S\'eminaire sure les Equations aux D\'eriv\'ees Partielles,
1992-1993, Ecole polytech., Palaiseau, 1993, Exp. No. XIII
\vol
\yr
\pages
\endref
\medskip

\ref
\key Gre
\by E. Grenier
\paper Semiclassical limit of the nonlinear Schr\"odinger equation in
small time
\jour Proc. Amer. Math. Soc.
\vol 126
\yr 1998
\issue 2
\pages 523-530
\endref
\medskip

\ref \key LZ \by F. Lin and P. Zhang
\paper Semiclassical limit of
the Gross-Pitaevskii equation in an exterior domain
\jour Arch.
Rational Mech. Anal.
\vol 179 \yr 2005
\issue
\pages 79--107
\endref

\medskip


\ref
\key RSS
\by I. Rodnianski, W. Schlag and  A. Soffer
\paper A symptotic stability of $N$-soliton states of NLS
\jour submitted
\vol
\yr
\pages
\endref

\medskip

\ref \key Sig \by I. M. Sigal \paper  "General characteristics of
non-linear dynamics,
\jour in "Spectral
 and
 Scattering Theory" (M. Ikawa, ed.), Marcel Dekker, Inc.
\yr 1994.
\pages
\endref

\medskip

\ref \key Sof \by  A. Soffer \paper Soliton Dynamics and Scattering
\jour International Congress of Mathematicians \vol III \yr 2006
\pages 459-471
\endref

\medskip

\ref \key T \by T. Tao
\paper A (concentration-)compact attractor
for high-dimensional non-linear Schrödinger equations \jour preprint
\vol \yr 2006 \pages
\endref

\medskip
\ref
\key
\by
\paper
\jour
\vol
\yr
\pages
\endref
\medskip

\endRefs

\newpage

\Refs
\widestnumber\key{CEL -Sog}

\ref
\key GV
\by J. Ginibre, G. Velo
\paper Generalized Strichartz Inequalities for the Wave Equation
\jour J. F. Analysis
\vol 133
\yr 1995
\pages 50--68
\endref
\medskip

\ref \key Gr \by J. M. Graf \paper Phase Space Analysis of the
Change transfer Model \jour HPA \vol 63\issue 2 \yr 1990 \pages
107--138
\endref
\medskip

\ref
\key JSS
\by J. L. Journ\'e, A. Soffer, C. D. Sogge
\paper Decay Estimates for the Schr\"odinger Operators
\jour Comm. Pure App. Math.
\vol 44
\yr 1991
\pages 573--604
\endref
\medskip

\ref
\key KT
\by M. Keel, T. Tao
\paper Endpoint Strichartz Estimates
\jour Amer. J. Math
\vol 120
\yr 1998
\pages 955--980
\endref
\medskip

\ref
\key NS
\by F. Nier, A. Soffer
\paper
\jour
\vol
\yr
\pages
\endref
\medskip



\ref
\key Pe
\by G. Perelman
\paper Some Results on the Scattering of Weakly Interacting Solitons
for Nonlinear Schr\"odinger Equations in ``Spectral theory, microload
analysis, singular manifolds''
\publ Akad. Verlag
\jour
\vol
\yr 1997
\pages 78--137
\endref
\medskip

\ref
\key RS
\by I. Rodnianski, W. Schlag
\paper Time Decay Solutions of Schr\"odinger Equations with Rough and
Time Dependent Potentials
\jour preprint (2001)
\publ
\vol
\yr
\pages
\endref
\medskip

\ref
\key St
\by R. S. Strichartz
\paper Restriction of Fourier Transforms to Quadratic Surfaces and
decay of Solutions of Wave Equations
\jour Duke Math J.
\vol 44\issue 3
\yr 1997
\pages 705--714
\endref
\medskip

\ref \key SW \by A. Soffer, M. Weinstein \paper Time Dependent
Resonance Theory \jour Geo. Func. Analysis (GAFA) \vol 8 \yr 1998
\pages 1086
\endref
\medskip



\ref
\key Ya 1
\by K. Yajima
\paper The $W^{k,p}$- continuity of Wave Operators for Schr\"odinger Operators
\jour J. Math. Soc., Japan
\publ
\vol 47\issue 3
\yr 1995
\pages 551--581
\endref
\medskip


\ref
\key Ya 2
\by K. Yajima
\paper A Multichannel Scattering Theory for some time dependent
\jour CMP
\vol 75\issue 2
\yr 1980
\pages 153--178
\endref
\medskip

\ref
\key Wu
\by U. W\"uller
\paper Geometric Methods in Scattering Theory of the Charge Transfer Model
\jour Duke Math J.
\vol 62
\issue 2
\yr 1991
\pages 273--313
\endref
\medskip

\ref
\key Z
\by L. Zielinski
\paper
\jour JFA
\vol 150
\issue 2
\yr 1997
\pages 453--470
\endref
\medskip


\endRefs
\enddocument

\end